# α-Q-fuzzy Subgroups

Muwafaq M. Salih and Delbrin H. Ahmed

Department of Mathmetics, College of Basic Education, Duhok University,

Duhok, Kurdistan Region – F.R. Iraq

ABSTRACT

In this paper, the notations of α-Q-fuzzy subset and α-Q-fuzzy subgroup are introduced, and a necessary properties related to these two concepts are proven. In the past part in this work, the effect of α-Q-fuzzy subgroup on the image and inverse-image under group anti-homomorphism are studied.

KEY WORDS: α-Q-fuzzy Subset, α-Q-fuzzy Subgroup, Anti-homomorphism, α-Q-fuzzy Abelian Subgroup and Cyclic α-Q-fuzzy Subgroup, AMS Subject Classifications (2010): 03E72, 08A72 and 20N25.

## 1. INTRODUCTION

In 1965, Zadeh first introduced the concept of fuzzy set. Then, a fuzzy subset developed widely in mathematics and was studied in several aspects in the mathematics, where in 1971 Rosenfeld used the concept of fuzzy sets in the fuzzy group theory. After that in 2006, Kim introduced the notation of intuitionistic Q-fuzzy semi-prime ideal in a semigroup. Solairaju and Nagarajan, in 2009, defined the notation of Q-fuzzy subgroup. Further, in 2013, Sharma [4] introduced new algebraic structure of α-fuzzy subgroup. Palaniappan and Muthuraj, in 2004, proved many results in homomorphism and anti-homomorphism in fuzzy subgroups.

The purpose of this work is to define notations of α-Q-fuzzy subset and α-Q-fuzzy subgroup, and we prove some elementary algebraic properties related to these two notations. Furthermore, the concept of anti-homomorphism in α-Q-fuzzy subgroup between two groups is defined and some properties related to an α-Q-fuzzy abelian and cyclic α-Q-fuzzy subgroups are studied.

## 2. PRELIMINARY

**Definition 2.1.** Zadeh (1965) - Let X be a non-empty set. A fuzzy subset θ of the set X is a mapping θ: X→[0, 1].

**Definition 2.2.** Rosenfeld (1971) - Let (G,∗) be a group. A fuzzy subset θ of G said to be a fuzzy subgroup of G, if the following conditions hold:
 i. θ(xy) ≥ min {θ(x), θ(y)} and
 ii. θ($x^{-1}$) ≥ θ(x), for all x and y in G.

**Definition 2.3.** Solairaju and Nagarajan (2006) - Let X and Q be non-empty sets. A Q-fuzzy subset θ of the set X is a mapping.
 i. θ: X×Q→[0, 1].

**Definition 2.4.** Solairaju and Nagarajan (2006) - Let (G,∗) be a group and Q be a non-empty set. A Q-fuzzy subset θ of G is said to be a Q-fuzzy subgroup of G if the following conditions are satisfied:
 i. θ(xy, q) ≥ min {θ(x, q), θ(y, q)} and
 ii. θ($x^{-1}$, q) ≥ θ(x, q), for all x and y in G and q in Q.

**Proposition 2.5.** Solairaju and Nagarajan (2006) - Let θ and σ be any two Q-fuzzy subsets of a non-empty set X. Then, for all x∈ X and q∈ Q, the followings are true:
 i. θ⊆σ ⇔ θ(x, q) ≤ σ(x, q),
 ii. θ=σ ⇔ θ(x, q) = σ(x, q),
 iii. (θ∪σ)(x, q) = max {θ(x, q), σ(x, q)},
 iv. (θ∩σ)(x, q) = min {θ(x, q), σ(x, q)}.

**Definition 2.6.** Solairaju and Nagarajan (2006) - If θ is a Q-fuzzy subgroup of a group G, then C[θ] is the complement of Q-fuzzy subgroup θ, and is defined by C[θ(x, q)] = 1-θ (x, q), for all x in G and q in Q.

**Definition 2.7.** Palaniappan and Muthuraj (2004) - If (G,∗) and (G′,o) are any two groups, then the function f: G→G′ is called a homomorphism if f (x∗y) = f(x) o f(y), for all x and y in G.

**Definition 2.8.** Palaniappan and Muthuraj (2004) If (G,∗) and (G′,o) are any two groups, then the function f: G→G′ is called an anti-homomorphism if f(x∗y) = f(y) o f(x), for all x and y in G.

**Definition 2.9.** Palaniappan and Muthuraj (2004) Let f: X→X′ be any mapping from a non-empty set X into a non-empty set X′, let θ be a fuzzy subset in X and σ be a fuzzy subset in X′ = f(X), defined by σ(x′) = sup θ(x), ∀ x∈ f⁻¹(x′), x∈ X and x′∈ X′. Then, θ is called an inverse image of σ under map f and is denoted by f⁻¹(σ).

**Definition 2.10.** Abdullah and Jeyaraman (2010) - If X = {x∈ G | θ(x) = θ(e)} where θ is fuzzy subgroups of a group G, then θ is fuzzy abelian subgroups of a group G.

**Definition 2.11.** Abdullah and Jeyaraman (2010) - Let G be a group and θ be a fuzzy subgroup of G. Then, θ is a cyclic fuzzy subgroup of G, if $θ_s$ is a cyclic subgroup for all s in [0, 1], and is defined as $θ_s$ = {x | θ(x) ≥ s, for x∈ G}.

**Definition 2.12.** Solairaju and Nagarajan (2006) - Let θ be a fuzzy subset of a group G. Let α∈ [0, 1]. Then, the fuzzy set $θ^α$ of G is called the α-fuzzy subset of G and is defined as $θ^α(x)$ = min {θ(x), α}, for all x∈ G.

**Definition 2.13.** Solairaju and Nagarajan (2006) - Let θ be a fuzzy subset of a group G. Let α∈ [0, 1]. Then, θ is called α-fuzzy subgroup of G and is denoted by $θ^α$, if for all x, y in G, the following conditions hold:
  i.  $θ^α(xy)$ ≥ min {$θ^α(x)$, $θ^α(y)$} and
  ii. $θ^α(x^{-1})$ ≥ $θ^α(x)$, ∀ x, y∈ G.

## 3. α-Q-FUZZY SUBSET

**Definition 3.1.** Let G and Q be any two non-empty sets and α∈ [0, 1]. Then, a mapping $θ^α$: G×Q→[0, 1] is called α-Q-fuzzy subset in G, with regard to fuzzy set θ, if $θ^α(x, q)$ = min {θ(x, q), α}, for x∈ G and q∈ Q.

**Remark 3.2.** Obviously, if α = 1, then $θ^α$ = θ, and if α = 0, then $θ^α$ = α.

**Proposition 3.3.** Let $θ^α$ and $σ^α$ be two α-Q-fuzzy subsets of G. Then, $(θ∩σ)^α$ = $θ^α ∩ σ^α$.

**Proof.** $(θ∩σ)^α(x, q)$ = min {(θ∩σ)(x, q), α}
= min {min {θ(x, q), σ(x, q)}, α}
= min {min {θ(x, q), α}, min σ(x, q), α}}
= min {$θ^α(x, q)$, $σ^α(x, q)$}
= $θ^α(x, q) ∩ σ^α(x, q)$
= $(θ^α ∩ σ^α)(x, q)$, for x∈ G and q∈ Q.

Hence, $(θ∩σ)^α$ = $θ^α ∩ σ^α$.

**Proposition 3.4.** Let f: G→G′ be a mapping. Let $θ^α$ be fuzzy subset of G and $σ^α$ be fuzzy subset of G′. Then;
  i.  $f(θ^α) = (f(θ))^α$,
  ii. $f^{-1}(σ^α) = (f^{-1}(σ))^α$.

**Proof.** (i) $f(θ^α(x, q))$ = sup {$θ^α(x, q)$ | f(x, q) = (x′, q)}
= sup {min {θ(x, q), α} | f(x, q) = (x′, q)}
= min {sup {θ(x, q) | f(x, q) = (x′, q)}, α}
= min {f(θ(x′, q)), α}
= $(f(θ)(x′, q))^α$, for all x∈ G, x′∈ G′ and q∈ Q.

Hence, $f(θ^α) = (f(θ))^α$.

(ii) $f^{-1}(σ^α(x, q)) = σ^α(f(x, q))$ = min {σ(f(x, q)), α}
= min {$f^{-1}(σ(x, q))$, α}
= $(f^{-1}(σ(x, q)))^α$, for all x∈ G and q∈ Q.

Hence, $f^{-1}(σ^α) = (f^{-1}(σ))^α$.

## 4. α-Q-FUZZY SUBGROUP

**Definition 4.1.** Let θ be a Q-fuzzy subgroup of a group G and α∈ [0, 1]. Then, $θ^α$ is called α-Q-fuzzy subgroup of G, if for all x, y∈ G and q∈ Q the following conditions hold:
  i.  $θ^α(xy, q)$ ≥ min {$θ^α(x, q)$, $θ^α(y, q)$} and
  ii. $θ^α(x^{-1}, q)$ ≥ $θ^α(x, q)$.

**Proposition 4.2.** Every Q-fuzzy subgroup of group G is α-Q-fuzzy subgroup of G.

**Proof.** Suppose θ be Q-fuzzy subgroup of G, and x, y be any two elements in a group G.

Then, $θ^α(xy, q)$ = min{θ(xy, q), α}
≥ min {min{θ(x, q), θ(y, q)}, α}
= min {min {θ(x, q), α}, min {θ(y, q), α}}
= min {$θ^α(x, q)$, $θ^α(y, q)$}

Therefore, $θ^α(xy, q)$ ≥ min {$θ^α(x, q)$, $θ^α(y, q)$}, and
$θ^α(x^{-1}, q)$ = min{θ($x^{-1}$, q), α} ≥ min{θ(x, q), α} = $θ^α(x, q)$

Therefore, $θ^α(x^{-1}, q)$ ≥ $θ^α(x, q)$.
Hence, θ is α-Q-fuzzy subgroup of G.

**Remark 4.3.** An α-Q-fuzzy subgroup of group G need not be Q-fuzzy subgroup of G.

**Note 4.4.** In all examples in this work, we consider Q = {q}.

**Example 4.5.** Conceder the four Klein group G= {e, a, b, c}, where $a^2 = b^2 = c^2$ = e and ab = ba = c, and a non-empty set Q = {q}. A Q-fuzzy subset θ of G is defined as:

$$(x,q) = \begin{cases} 0.2 & \text{if } x=e \\ 0.4 & \text{if } x=a \text{ or } b \\ 0.3 & \text{if } x=c \end{cases}$$

It is clear, the Q-fuzzy subset θ is not a Q-fuzzy Subgroup of G since first part of definition fails to hold as θ(c, q) < min {θ(a, q), θ(b, q)}.

However, it can be shown that θ is α- Q- fuzzy subgroups of group G. This, if we take α = 0.09, it can be observed, θ(x, q) > α, for all elements in the group G and q in Q. This implies that, θ(x, q) = min {θ(x, q), α} = α, ∀x∈ G and q∈ Q, which satisfies the first part of the definition of α-Q-fuzzy subgroups, $θ^α(xy, q)$ ≥ min {θ(x, q), θ(x, q)}. Furthermore, the second part of the definition, since $a^{-1}$= a, $b^{-1}$= b and $c^{-1}$=c, this satisfies $θ^α(x^{-1}, q)$ ≥ $θ^α(x, q)$. Hence, $θ^α$ is α- Q- fuzzy subgroups of group G, and denoted by $θ^α$.

**Proposition 4.6.** If $θ^α$ is α-Q-fuzzy subgroup of a group G, then
  i.  $θ^α(e, q)$ ≥ $θ^α(x, q)$ for all x∈ G and q∈ Q, where e is the identity element of a group G.
  ii. A set K = {x∈ G | $θ^α(x, q)$ = $θ^α(e, q)$, for q∈ Q} is an α-Q-fuzzy subgroup of a group G.

**Proof.** (i) $θ^α(e, q)$ = $θ^α(xx^{-1}, q)$ ≥ min {$θ^α(x, q)$, $θ^α(x^{-1}, q)$}
= min {$θ^α(x, q)$, $θ^α(x, q)$}
= $θ^α(x, q)$.

Hence, $θ^α(e, q)$ ≥ $θ^α(x, q)$ ∀ x∈ G and q∈ Q.

(ii) The set K ≠ ∅ as at least there exists e∈ K. Let x, y∈ K and q∈ Q. Then, $θ^α(x, q)$ = $θ^α(y, q)$ = $θ^α(e, q)$.

$θ^α(xy^{-1}, q)$ ≥ min{$θ^α(x, q)$, $θ^α(y^{-1}, q)$}
= min {$θ^α(x, q)$, $θ^α(y, q)$}

$\qquad = \min \{\theta^\alpha(e, q), \theta^\alpha(e, q)\}$
$\qquad = \theta^\alpha(e, q)$
So, $\theta^\alpha(xy^{-1}, q) \geq \theta^\alpha(e, q)$.
By part (i), we can show that $\theta^\alpha(e, q) \geq \theta^\alpha(xy^{-1}, q)$
Therefore, $\theta^\alpha(xy^{-1}, q) = \theta^\alpha(e, q)$.
Hence, K is an α-Q-fuzzy subgroup of the group G.

**Proposition 4.7.** Let $\theta^\alpha$ be α-Q-fuzzy subgroups of a group G. If $\theta^\alpha(xy^{-1}, q) = \theta^\alpha(e, q)$, then $\theta^\alpha(x, q) = \theta^\alpha(y^{-1}, q)$, for all x and y in G and q in Q.

**Proof.** Let $\theta^\alpha$ be α-Q-fuzzy subgroup of a group G, and $\theta^\alpha(xy^{-1}, q) = \theta^\alpha(e, q)$, for all x, y ∈ G and q ∈ Q.
Then, $\theta^\alpha(x, q) = \theta^\alpha(x(y^{-1}y), q) = \theta^\alpha((xy^{-1})y, q)$
$\qquad \geq \min \{\theta^\alpha(xy^{-1}, q), \theta^\alpha(y, q)\}$
$\qquad = \min \{\theta^\alpha(e, q), \theta^\alpha(y, q)\}$
$\qquad = \theta^\alpha(y, q)$.
Therefore, $\theta^\alpha(x, q) \geq \theta^\alpha(y, q)$, and
$\theta^\alpha(y, q) = \theta^\alpha(y^{-1}, q) = \theta^\alpha((x^{-1}x)y^{-1}, q) = \theta^\alpha(x^{-1}(xy^{-1}), q)$
$\qquad \geq \min \{\theta^\alpha(x^{-1}, q), \theta^\alpha(xy^{-1}, q)\}$
$\qquad = \min \{\theta^\alpha(x, q), \theta^\alpha(e, q)\}$
$\qquad = \theta^\alpha(x, q)$.
Therefore, $\theta^\alpha(y, q) \geq \theta^\alpha(x, q)$.
Hence, $\theta^\alpha(x, q) = \theta^\alpha(y, q)$.

**Proposition 4.8.** If $\theta^\alpha$ and $\sigma^\alpha$ be two α-Q-fuzzy subgroups of a group G, then $(\theta \cap \sigma)^\alpha$ is also an α-Q-fuzzy subgroup of a group G.

**Proof.** Let x, y be any two elements in G and q in Q. Then,
$(\theta \cap \sigma)^\alpha(xy, q) = (\theta^\alpha \cap \sigma^\alpha)(xy, q)$
$\qquad = \min \{\theta^\alpha(xy, q), \sigma^\alpha(xy, q)\}$
$\qquad \geq \min \{\min\{\theta^\alpha(x, q), \theta^\alpha(y, q)\}, \min\{\sigma^\alpha(x, q), \sigma^\alpha(y, q)\}\}$
$\qquad = \min \{\min\{\theta^\alpha(x, q), \sigma^\alpha(x, q)\}, \min \{\theta^\alpha(y, q), \sigma^\alpha(y, q)\}\}$
$\qquad = \min \{(\theta \cap \sigma)^\alpha(x, q), (\theta \cap \sigma)^\alpha(y, q)\}$
Therefore, $(\theta \cap \sigma)^\alpha(xy, q) \geq \min \{(\theta \cap \sigma)^\alpha(x, q), (\theta \cap \sigma)^\alpha(y, q)\}$, and
$(\theta \cap \sigma)^\alpha(x^{-1}, q) = (\theta^\alpha \cap \sigma^\alpha)(x^{-1}, q)$
$\qquad = \min \{\theta^\alpha(x^{-1}, q), \sigma^\alpha(x^{-1}, q)\}$
$\qquad \geq \min \{\theta^\alpha(x, q), \sigma^\alpha(x, q)\}$
$\qquad = (\theta^\alpha \cap \sigma^\alpha)(x, q)$.
Therefore, $(\theta \cap \sigma)^\alpha(x^{-1}, q) \geq (\theta^\alpha \cap \sigma^\alpha)(x, q)$.
Hence, $(\theta \cap \sigma)^\alpha$ is α-Q-fuzzy subgroup of a group G.

**Remark 4.9.** If $\theta^\alpha$ and $\sigma^\alpha$ be two α-Q-fuzzy subgroup of a group G, then $(\theta \cup \sigma)^\alpha$ need not be an α-Q-fuzzy subgroup of a group G.

**Example 4.10.** Consider G be a group of integer under addition, and θ, σ and π be Q-fuzzy sets of the group G defined as:

$(x,q)=\begin{cases} 0.4 & \text{if } x \in 3Z \\ 0 & \text{otherwise} \end{cases}, (x,q)=\begin{cases} 0.2 & \text{if } x \in Z_e \\ 0.1 & \text{otherwise} \end{cases}$

and $p(x,q)=\begin{cases} 1 & \text{if } x \in Z_e \\ 0 & \text{otherwise} \end{cases}$

It can be shown that, the θ, σ and π are α-Q-fuzzy subgroups of a group G, and denoted by $\theta^\alpha$, $\sigma^\alpha$ and $\pi^\alpha$, respectively.

However, the union of two α-Q-fuzzy subgroups θ and σ is not α-Q-fuzzy subgroups, where $(\theta \cup \sigma)(x, q) = \max\{\theta(x, q), \sigma(x, q)\}$. Therefore,

$(q \cup s)(x,q)=\begin{cases} 0.4 & \text{if } x \in 3Z \\ 0.2 & \text{if } x \in Z_e \\ 0.1 & \text{otherwise} \end{cases}$

Now, if we take x = 3 and y = 2, then $(\theta \cup \sigma)(x, q) = \max \{\theta(x, q), \sigma(x, q)\} = 0.4$, $(\theta \cup \sigma)(y, q) = \max \{\theta(y, q), \sigma(y, q)\} = 0.2$, and $(\theta \cup \sigma)(x-y, q) = (\theta \cup \sigma)(3-2, q) = (\theta \cup \sigma)(1, q) = 0.1$. Therefore, the first part of definition fails to hold, as $(\theta \cup \sigma)(x-y, q) = 0.1 < \min\{(\theta \cup \sigma)(x, q), (\theta \cup \sigma)(y, q)\} = 0.2$. Hence, $(\theta \cup \sigma)$ is not α-Q-fuzzy subgroup of the group G.

On the other hand, the union of Q-fuzzy subgroups σ and π is $(\sigma \cup \pi)(x,q)=\begin{cases} 1 & \text{if } x \in Z_e \\ 0.1 & \text{otherwise} \end{cases}$

It can be observed that, $(\sigma \cup \pi)$ is α-Q-fuzzy subgroup of the group G, and denoted by $(\sigma \cup \pi)^\alpha$.

To conclude above, the union of two α-Q-fuzzy subgroups of the group G need not be an α-Q-fuzzy subgroup of G.

**Proposition 4.11.** If $\theta^\alpha$ is an α-Q-fuzzy subgroup of a group G, then $C(\theta^\alpha)$ is an α-Q-fuzzy subgroup of the group G.

**Proof.** $C[\theta^\alpha(xy, q)] = 1-\theta^\alpha(xy, q)$
$\qquad \leq 1-\min \{\theta^\alpha(x, q), \theta^\alpha(y, q)\}$
$\qquad = \max \{1-\theta^\alpha(x, q), 1-\theta^\alpha(y, q)\}$
$\qquad = \max \{C[\theta^\alpha(x, q)], C[\theta^\alpha(y^{-1}, q)]\}$
Therefore, $C[\theta^\alpha(xy, q)] \leq \max \{C[\theta^\alpha(x, q)], C[\theta^\alpha(y^{-1}, q)]\}$, ∀ x ∈ G and q ∈ Q, and
$C[\theta^\alpha(x^{-1}, q)] = 1-\theta^\alpha(x^{-1}, q)$
$\qquad \geq 1-\theta^\alpha(x, q)$
$\qquad = C[\theta^\alpha(x, q)]$
Therefore, $C[\theta^\alpha(x^{-1}, q)] \geq = C[\theta^\alpha(x, q)]$
Hence, $C(\theta^\alpha)$ is an α-Q-fuzzy subgroup of the group G.

**Proposition 4.12.** Let θ is Q-fuzzy subgroup of a group G. Then, $\theta^\alpha$ is α-Q-fuzzy subgroups of a group G if and only if $\theta^\alpha(xy^{-1}, q) \geq \min \{\theta^\alpha(x, q), \theta^\alpha(y, q)\}$ for all x, y in G and q in Q.

**Proof.** Suppose $\theta^\alpha$ is α-Q-fuzzy subgroup of a group G.
Then, $\theta^\alpha(xy^{-1}, q) \geq \min \{\theta^\alpha(x, q), \theta^\alpha(y^{-1}, q)\}$
$\qquad \geq \min \{\theta^\alpha(x, q), \theta^\alpha(y, q)\}$.
Therefore, $\theta^\alpha(xy^{-1}, q) \geq \min \{\theta^\alpha(x, q), \theta^\alpha(y, q)\}$, for all x, y in G and q in Q.
Conversely, suppose $\theta^\alpha(xy^{-1}, q) \geq \min \{\theta^\alpha(x, q), \theta^\alpha(y, q)\}$, for all x, y in G and q in Q.
Then, $\theta^\alpha(xy, q) = \theta^\alpha(x(y^{-1})^{-1}, q)$
$\qquad \geq \min \{\theta^\alpha(x, q), \theta^\alpha(y^{-1}, q)\}$
$\qquad \geq \min \{\theta^\alpha(x, q), \theta^\alpha(y, q)\}$.
Therefore, $\theta^\alpha(xy, q) \geq \min \{\theta^\alpha(x, q), \theta^\alpha(y, q)\}$, for all x, y in G and q in Q, and
$\theta^\alpha(x^{-1}, q) = \theta^\alpha(x^{-1}e, q) \geq \min \{\theta^\alpha(x, q), \theta^\alpha(e, q)\} = \theta^\alpha(x, q)$.
Therefore, $\theta^\alpha(x^{-1}, q) \geq \theta^\alpha(x, q)$, for all x in G and q in Q, where e is the identity element of G.
Hence, $\theta^\alpha$ is α-Q-fuzzy subgroup of a group G.

**Definition 4.13.** Let $\theta^\alpha$ and $\sigma^\alpha$ be any two α-Q-fuzzy subgroups of groups G and G′, respectively. Then, $\theta^\alpha \times \sigma^\alpha$ said to be product of $\theta^\alpha$ and $\sigma^\alpha$ and is defined as $\theta^\alpha \times \sigma^\alpha ((x, x'), q) = \min \{\theta^\alpha(x, q), \sigma^\alpha(x', q)\}$, $\forall x \in G$, $x' \in G'$ and $q \in Q$.

**Proposition 4.14.** If $\theta^\alpha$ and $\sigma^\alpha$ be two α-Q-fuzzy subgroups of groups G and G′ respectively, then $\theta^\alpha \times \sigma^\alpha$ is α-Q-fuzzy subgroup of a group G×G′.

**Proof.** If $x, y \in G$, and $x', y' \in G'$, then $(x, x'), (y, y') \in G \times G'$,

$\theta^\alpha \times \sigma^\alpha((x, x')(y, y'), q) = \theta^\alpha \times \sigma^\alpha([(xy),(x'y')], q)$
$= \min \{\theta^\alpha(xy, q), \sigma^\alpha(x'y', q)\}$
$\geq \min\{\min\{\theta^\alpha(x,q), \theta^\alpha(y,q)\}, \min\{\sigma^\alpha(x', q), \sigma^\alpha(y', q)\}\}$
$= \min \{\min\{\theta^\alpha(x, q), \sigma^\alpha(x', q)\ \theta^\alpha(y, q)\}, \min\{\theta^\alpha(y, q), \sigma^\alpha(y', q)\}\}$
$= \min \{\theta^\alpha \times \sigma^\alpha ((xy), q), \theta^\alpha \times \sigma^\alpha ((xy), q)\}$

Therefore, $\theta^\alpha \times \sigma^\alpha ((x, x')(y, y'), q) \geq \min \{\theta^\alpha \times \sigma^\alpha((xy), q), \theta^\alpha \times \sigma^\alpha((xy), q)\}$, and

$\theta^\alpha \times \sigma^\alpha((x, x')^{-1}, q) = \theta^\alpha \times \sigma^\alpha((x^{-1}, x'^{-1}), q)$
$= \min \{\theta^\alpha(x^{-1}, q), \sigma^\alpha(x'^{-1}, q)\}$
$\geq \min \{\theta^\alpha(x, q), \sigma^\alpha(x', q)\}$
$= \theta^\alpha \times \sigma^\alpha((x, x'), q)$

Therefore, $\theta^\alpha \times \sigma^\alpha((x, x')^{-1}, q) \geq \theta^\alpha \times \sigma^\alpha((x, x'), q)$

Hence, $\theta^\alpha \times \sigma^\alpha$ is α-Q-fuzzy subgroup of a group G×G′.

**Proposition 4.15.** If $\theta^\alpha$, $\sigma^\alpha$ and $\theta^\alpha \times \sigma^\alpha$ are α-Q-fuzzy subgroups of groups G, G′ and G×G′, respectively, then the following statements are true:
(i) If $\theta^\alpha(e, q) \geq \sigma^\alpha(x, q)$ $\forall x \in G$ and $q \in Q$, then $\sigma^\alpha$ is an α-Q-fuzzy subgroup of G′.
(ii) If $\sigma^\alpha(e', q) \geq \theta^\alpha(x, q)$ $\forall x \in G$ and $q \in Q$, then $\theta^\alpha$ is an α-Q-fuzzy subgroups of G.

Where e and e′ are the identity elements of groups G and G′, respectively.

**Proof.** Let $\theta^\alpha \times \sigma^\alpha$ be an α-Q-fuzzy subgroup of G×G′ with $x, y \in G$ and $q \in Q$.

(i) Since $\theta^\alpha(e, q) \geq \sigma^\alpha(x, q)$, $\forall x \in G$ and $q \in Q$.

Then, $\sigma^\alpha(xy^{-1}, q) = \min \{\theta^\alpha(ee, q), \sigma^\alpha(xy^{-1}, q)\}$
$= \theta^\alpha \times \sigma^\alpha([(ee), (xy^{-1})], q)$
$= \theta^\alpha \times \sigma^\alpha([(e, x), (e, y^{-1})], q)$
$\geq \min \{\theta^\alpha \times \sigma^\alpha((e, x), q), \theta^\alpha \times \sigma^\alpha((e, y^{-1}), q)\}$
$= \min \{\min \{\theta^\alpha(e, q), \sigma^\alpha(x, q)\}, \min \{\theta^\alpha(e, q), \sigma^\alpha(y^{-1}, q)\}\}$
$= \min \{\sigma^\alpha(x, q)\}, \sigma^\alpha(y^{-1}, q)\}$
$\geq \min \{\sigma^\alpha(x, q)\}, \sigma^\alpha(y, q)\}$

Therefore, $\sigma^\alpha(xy^{-1}, q) \geq \min \{\sigma^\alpha(x, q), \sigma^\alpha(y, q)\}$, $\forall x, y \in G$ and $q \in Q$.

Hence, $\sigma^\alpha$ is an α-Q-fuzzy subgroup of G′.

(ii) Since $\sigma^\alpha(e', q) \geq \theta^\alpha(x, q)$, e′ is the identity element of G′, $\forall x \in G$ and $q \in Q$.

Then, $\theta^\alpha(xy^{-1}, q) = \min \{\theta^\alpha(xy^{-1}, q), \sigma^\alpha((e'e'), q)\}$
$= \theta^\alpha \times \sigma^\alpha([(xy^{-1}), (e'e')], q)$
$= \theta^\alpha \times \sigma^\alpha([(x, e'), (y^{-1}, e')], q)$
$\geq \min \{\theta^\alpha \times \sigma^\alpha((x, e'), q), \theta^\alpha \times \sigma^\alpha((y^{-1}, e'), q)\}$
$= \min \{\min \{\theta^\alpha(x, q), \sigma^\alpha(e', q)\}, \min \{\theta^\alpha(y^{-1},q), \sigma^\alpha(e', q)\}\}$
$= \min \{\theta^\alpha(x, q)\}, \theta^\alpha(y^{-1}, q)\}$
$\geq \min \{\theta^\alpha(x, q)\}, \theta^\alpha(y, q)\}$

Therefore, $\theta^\alpha(xy^{-1}, q) \geq \min \{\theta^\alpha(x, q), \theta^\alpha(y, q)\}$, $\forall x, y \in G$ and $q \in Q$.

Hence, $\theta^\alpha$ is an α-Q-fuzzy subgroup of G.

**Remark 4.16.** If $\theta^\alpha$, $\sigma^\alpha$ and $\theta^\alpha \times \sigma^\alpha$ are α-Q-fuzzy subgroups of G, G′ and G×G′, respectively, then either $\theta^\alpha$ is a α-Q-fuzzy subgroup of G or $\sigma^\alpha$ is a α-Q-fuzzy subgroup of G′.

**Proof.** The proof is clear from proposition 4.15.

## 5. α-Q-FUZZY SUBGROUPS UNDER ANTI-HOMOMORPHISM

**Definition 5.1.** Let f: G→G′ be any mapping from a group G into a group G′, let $\theta^\alpha$ be an α-Q-fuzzy subgroup in G and $\sigma^\alpha$ be an α-Q-fuzzy subgroup in G′ = f(G), defined by $\sigma^\alpha(x', q) = \sup \theta^\alpha(x, q)$, $\forall x \in f^{-1}(x')$, $x \in X$ and $x' \in X'$. Then, $\theta^\alpha$ is called an inverse image of $\sigma^\alpha$ under map f and is denoted by $f^{-1}(\sigma^\alpha)$.

**Proposition 5.2.** Let f: G→G′ be an anti-homomorphism from a group G into the group G′, and Q be a non-empty set. If $\theta^\alpha$ is an α-Q-fuzzy subgroup of G, then the anti-homomorphic image $f(\theta^\alpha)$ is an α-Q-fuzzy subgroup of G′.

**Proof.** Let $\theta^\alpha$ be an α-Q-fuzzy subgroup of G and its image $f(\theta^\alpha)$ be in G′.

Then for all f(x), f(y) in G, q in Q,
$f(\theta^\alpha(f(x)f(y), q)) = f(\theta^\alpha(f(yx), q))$
$= \theta^\alpha(yx, q)$
$\geq \min \{\theta^\alpha(x, q), \theta^\alpha(y, q)\}$
$= \min \{f(\theta^\alpha(f(x), q)), f(\theta^\alpha f(y), q))\}$

Therefore, $f(\theta^\alpha(f(x)f(y), q)) \geq \min \{f(\theta^\alpha(f(x), q)), f(\theta^\alpha f(y), q))\}$, and, for all f(x) in G′ and q in Q.
$f(\theta^\alpha(f(x)^{-1}, q)) = f(\theta^\alpha(f(x^{-1}), q))$
$= \theta^\alpha(x^{-1}, q)$
$\geq \theta^\alpha(x^{-1}, q)$
$= f(\theta^\alpha(f(x), q))$

Therefore, $f(\theta^\alpha(f(x)^{-1}, q)) \geq f(\theta^\alpha(f(x), q))$

Hence, $f(\theta^\alpha)$ is an α-Q-fuzzy subgroup of G′.

**Remark 5.3.** Let f: G→G′ be a homomorphism from a group G into the group G′, and Q be a non-empty set. If $\theta^\alpha$ is an α-Q-fuzzy subgroup of G, then the homomorphic image $f(\theta^\alpha)$ is an α-Q-fuzzy subgroup of G′.

**Proof.** Follows from proposition 5.2.

**Proposition 5.4.** Let f: G→G′ be an anti-homomorphism from a group G into the group G′, and Q be a nonempty set. If $\sigma^\alpha$ is an α-Q-fuzzy subgroup of f(G) = G′, then the anti-homomorphic inverse image $f^{-1}(\sigma^\alpha)$ is an α-Q-fuzzy subgroup of G.

**Proof.** Let $\sigma^\alpha$ be an α-Q-fuzzy subgroup of G′ and x, y in G, q in Q.

Then, $f^{-1}(\sigma^\alpha(xy, q)) = \sigma^\alpha[f(xy), q]$
$= \sigma^\alpha[(f(y)f(x), q]$
$\geq \min \{\sigma^\alpha[(f(y), q], \sigma^\alpha[(f(x), q]\}$
$= \min \{\sigma^\alpha[(f(x), q], \sigma^\alpha[(f(y), q]\}$
$= \min \{f^{-1}(\sigma^\alpha(x, q), f^{-1}(\sigma^\alpha(y, q)\}$

Therefore, $f^{-1}(\sigma^\alpha(xy, q)) \geq \min \{f^{-1}(\sigma^\alpha(y, q), f^{-1}(\sigma^\alpha(x, q)\}$, and, for all x in G and q in Q,
$f^{-1}(\sigma^\alpha(x^{-1}, q) = \sigma^\alpha[f(x^{-1}), q]$

$$\geq \sigma^\alpha[f(x), q]$$
$$= f^{-1}(\sigma^\alpha(x, q)).$$

Therefore, $f^{-1}(\sigma^\alpha(x^{-1}, q)) \geq f^{-1}(\sigma^\alpha(x, q))$.

Hence, $f^{-1}(\sigma^\alpha)$ is an α-Q-fuzzy subgroup of G.

**Remark 5.5.** Let f: G→G' be a homomorphism from a group G into the group G', and Q be a non-empty set. If $\sigma^\alpha$ is an α-Q-fuzzy subgroup of f(G) = G', then the homomorphic inverse image $f^{-1}(\sigma^\alpha)$ is an α-Q-fuzzy subgroup of G.

**Proof.** Follows from proposition 5.4.

**Definition 5.6.** If X = {x∈G | $\theta^\alpha(x, q) = \theta^\alpha(e, q)$, for q∈Q} where $\theta^\alpha$ is α-Q-fuzzy subgroup of a group G, then $\theta^\alpha$ is α-Q-fuzzy abelian subgroup of a group G.

**Proposition 5.7.** An Anti-homomorphism image of an α-Q-fuzzy abelian subgroups is an α-Q-fuzzy abelian subgroups.

**Proof.** Assume, f: G→G' is an anti-homomorphism and $\sigma^\alpha$ is an α-Q-fuzzy subgroup of G'.

Since $\theta^\alpha$ is α-Q-fuzzy abelian subgroup of the group G, implies that X ={x∈G | $\theta^\alpha(x, q) = \theta^\alpha(e_1, q)$, for q∈Q} is α-Q-fuzzy abelian subgroup of the group G, where $e_1$ the identity element of G.

Now consider, X' ={x'∈G' | $\sigma^\alpha(x, q) = \sigma^\alpha(e_2, q)$, for q∈Q} is α-Q-fuzzy subgroup of the group G', where $e_2$ is the identity element in G'. Given that $\sigma^\alpha$ is an α-Q-fuzzy subgroup of G'.

If xx'∈ X' and q∈Q, then $\sigma^\alpha(xx', q) = \sigma^\alpha(e_2, q) \Rightarrow \sup\theta^\alpha(\iota, q)$ for $\iota \in f^{-1}(xx', q) = \sup\theta^\alpha(\iota, q)$ for $\iota \in f^{-1}(e_2, q)$
$$\Rightarrow \theta^\alpha(xx', q) = \theta^\alpha(e_1, q)$$

Therefore, xx'∈ X and X is abelian.

Hence, $\theta^\alpha(xx', q) = \theta^\alpha(x'x, q) \Rightarrow \sup \theta^\alpha(\iota, q)$ for $\iota \in f^{-1}(xx', q)$
$$= \sup\theta^\alpha(\iota, q) \text{ for } \iota \in f^{-1}(x'x, q)$$
$$\Rightarrow \sigma^\alpha(xx', q) = \sigma^\alpha(x'x, q)$$
$$\Rightarrow \sigma^\alpha(e_2, q) = \sigma^\alpha(x'x, q)$$
$$\Rightarrow x'x \in X'.$$

Therefore, xx' = x'x in X'.

Hence, X' is abelian, and thus $\sigma^\alpha$ is α-Q-fuzzy abelian subgroup of the group G'.

**Proposition 5.8.** An anti-homomorphism inverse-image of an α-Q-fuzzy abelian subgroup is an α-Q-fuzzy abelian subgroup.

**Proof.** Assume $\theta^\alpha$ is α-Q-fuzzy subgroup of the group G and $\sigma^\alpha$ is α-Q-fuzzy abelian subgroup of the group G'. Then, X'={y∈G' | $\sigma^\alpha(y, q) = \sigma^\alpha(e_2, q)$, for q∈Q} is α-Q-fuzzy abelian subgroup of the group G', where $e_2$ is the identity element of G'.

Now consider, X= {x∈G | $\theta^\alpha(x, q) = \theta^\alpha(e_1, q)$, for q∈Q} is α-Q-fuzzy subgroup of the group G, where $e_1$ is the identity element of G.

If xx'∈ X and q∈Q, then, $\theta^\alpha(xx', q) = \theta^\alpha(e_1, q)$
$$\Rightarrow \sigma^\alpha[f(xx', q)] = \sigma^\alpha[f(e_1, q)]$$
$$\Rightarrow \sigma^\alpha[f(xx', q)] = \sigma^\alpha[f(e_2, q)]$$
$$\Rightarrow \sigma^\alpha[f(x, q)f(x', q)] = \sigma^\alpha[f(e_2, q)]$$

Therefore, f(x, q)f(x', q)∈ X', and since X' is abelian, this implies that f(x, q)f(x', q) = f(x', q)f(x, q).

Then, $\sigma^\alpha[f(x, q)f(x', q)] = \sigma^\alpha[f(x', q)f(x, q)]$
$$\Rightarrow \sigma^\alpha[f(xx', q)] = \sigma^\alpha[f(x'x, q)]$$

$$\Rightarrow \theta^\alpha(xx', q) = \theta^\alpha(x'x, q)$$
$$\Rightarrow \theta^\alpha(e_1, q) = \theta^\alpha(x'x, q)$$
$$\Rightarrow x'x \in X.$$

Therefore, xx' = x'x in X.

Hence, X is abelian, and thus $\theta^\alpha$ is α-Q-fuzzy abelian subgroup of the group G.

**Definition 5.9.** Let G be a group and $\theta^\alpha$ be an α-Q-fuzzy subgroup of G. Then, $\theta^\alpha$ is a cyclic α-Q-fuzzy subgroup of G, if $\theta_s$ is a cyclic subgroup for all s in [0, 1], where $\theta_c^\alpha$ = {(x, q) | $\theta^\alpha(x, q) \geq c$, for x∈G and q∈Q}.

**Proposition 5.10.** An anti-homomorphism image of a cyclic α-Q-fuzzy subgroup is a cyclic α-Q-fuzzy subgroup.

**Proof.** Assume $\sigma^\alpha$ is α-Q-fuzzy subgroup of the group G'. Our claim is to prove $\sigma^\alpha$ is a cyclic α-Q-fuzzy subgroup of G. Consider, f: G→G' is an anti-homomorphism and $\theta^\alpha$ is a cyclic α-Q-fuzzy subgroup of G.

Obviously, for any c∈[0, 1], a $\theta_c^\alpha$ is a cyclic subgroup of G, and f($\theta_c^\alpha$) ⊆ $\sigma_c^\alpha$. Therefore, $\sigma_c^\alpha$ is cyclic, and hence $\sigma^\alpha$ is cyclic α-Q-fuzzy subgroup of G'.

**Proposition 5.11.** An anti-homomorphism inverse image of a cyclic α-Q-fuzzy subgroup is a cyclic α-Q-fuzzy subgroup.

**Proof.** Assume $\theta^\alpha$ is α-Q-fuzzy subgroup of the group G. Our claim is to prove $\theta^\alpha$ is a cyclic α-Q-fuzzy subgroup of G. Consider, f: G→G' is an anti-homomorphism and $\sigma^\alpha$ is a cyclic α-Q-fuzzy subgroup of G'.

Now, for any c∈[0, 1], $\theta_c^\alpha$ = {(x, q) | $\theta^\alpha(x, q) \geq c$, for x∈G and q∈Q}. Since $\sigma^\alpha$ is a cyclic α-Q-fuzzy subgroup of G', then $\sigma_c^\alpha$ is a cyclic subgroup of G'. So, obviously $f^{-1}(\sigma_c^\alpha) \subseteq \theta_c^\alpha$. Therefore, $\theta_c^\alpha$ is cyclic, and hence $\theta^\alpha$ is cyclic α-Q-fuzzy subgroup of G.

## 6. CONCLUSION

In this paper, the concepts of α-Q-fuzzy subset and α-Q-fuzzy subgroup have been defined and related properties are proven. Furthermore, a group anti-homomorphism has been in terms of image and inverse image and its effect to α-Q-fuzzy abelian subgroup and cyclic α-Q-fuzzy subgroup are studied. In the following work, we may introduce the concept of α-Q-fuzzy normal subgroup and its effects on image and inverse image under homomorphism and anti-homomorphism could be studied.